\documentclass[12pt,noamsfonts]{amsart}

\usepackage{amsfonts}
\usepackage{amssymb}
\usepackage{amscd}

\newcommand{\marginlabel}[1]%
  {\mbox{}\marginpar{\raggedleft\hspace{0pt}\bfseries\sf#1}}
\def\NN{{\mathbb N}}

\def\CC{{\mathbb C}}
\def\AA{{\mathbb A}}
\def\RR{{\mathbb R}}
\def\QQ{{\mathbb Q}}

\def\cI{\mathcal{I}}

\def\cF{\mathcal{F}}

\def\cO{\mathcal{O}}

\def\ua{\underline{\mathbf{a}}}
\def\ub{\underline{\mathbf{b}}}
\def\uc{\underline{\mathbf{c}}}
\def\ur{\underline{\mathbf{r}}}
\def\um{\underline{\mathbf{m}}}

\newtheorem{lemma}{Lemma}[section]
\newtheorem{theorem}[lemma]{Theorem}
\newtheorem{corollary}[lemma]{Corollary}
\newtheorem{proposition}[lemma]{Proposition}

\theoremstyle{definition}

\newtheorem{remark}[lemma]{Remark}

\theoremstyle{remark}
\newtheorem*{remark*}{Remark}

\begin{document}

\title{The multiplier ideals of a sum of ideals}

\author[M. Musta\c{t}\v{a}]{Mircea~Musta\c{t}\v{a}}
\address{Department of Mathematics, University of California,
Berkeley, CA, 94720 and Institute of Mathematics of the Romanian
Academy}
\email{{\tt mustata@math.berkeley.edu}}

\thanks{\mbox {  } \mbox {  } 2000
\newline \mbox {  } \mbox {  }
}

\begin{abstract}
We prove that if $\ua$, $\ub\subseteq\cO_X$ are nonzero sheaves of ideals
on a complex smooth variety $X$, then for every $\gamma\in\QQ_+$
we have the following relation between the multiplier ideals of
$\ua$, $\ub$ and $\ua+\ub$:
$$\cI\left(X,\gamma\cdot(\ua+\ub)\right)\subseteq\sum_{\alpha+\beta=\gamma}
\cI(X,\alpha\cdot\ua)\cdot\cI(X,\beta\cdot\ub).$$
A similar formula holds for the asymptotic mutiplier ideals of the sum of two 
graded systems of ideals.

We use this result to approximate at a given point arbitrary multiplier ideals
by multiplier ideals associated to zero dimensional ideals. This is applied
to compare the multiplier ideals associated to a scheme in different
embeddings.
\end{abstract}

\maketitle

\bigskip

\section*{Introduction}

Let $X$ be a smooth complex variety. To a nonzero 
quasi-coherent sheaf of ideals
$\ua$ on $X$ one can associate a sequence of ideals called the multiplier
 ideals of $\ua$, which depend on a rational parameter. The behaviour of
 these ideals encodes in a subtle way the properties of the singularities
of $V(\ua)$.
Introduced first in the analytic context in the work of Demailly, Nadel, Siu
and others, multiplier ideals have recently found surprising applications
in algebraic geometry (see \cite{ein}, \cite{siu}, \cite{kawamata1},
\cite{kawamata2}, \cite{el}, \cite{els}).

Here is the definition. Suppose that $f\,:\,X'\longrightarrow X$ is a log 
resolution of $(X, V(\ua))$ i.e. $f$ is proper and birational, $X'$ is smooth,
and $f^{-1}V(\ua)=D$ is a divisor with simple normal crossings.
If $K_{X'/X}$ is the relative canonical divisor of $f$, the multiplier ideal of
$\ua$ with coefficient $\alpha\in\QQ_+$ is
$$\cI(X, \alpha\cdot\ua)=f_*\cO_{X'}(K_{X'/X}-[\alpha D]).$$
Here $[\cdot]$ denotes the integral part function.

In general one can expect that algebraic properties of the multiplier ideals
are related to the behaviour of linear systems and singularities of algebraic
varieties. For example, in \cite{del} is proved the following subadditivity
relation:
 $$\cI(X, \alpha\cdot(\ua\cdot\ub))\subseteq\cI(X,\alpha\cdot\ua)\cdot
\cI(X, \alpha\cdot\ub).$$ 
This is applied in \cite{del} to prove a theorem of Fujita on the volume
of a big divisor and in \cite{els} to show a uniform behaviour of the
symbolic powers of an ideal.

\smallskip

The main result of this paper is an analogous formula for the sum 
of two ideals.

\begin{theorem}\label{one_ideal}
If $X$ is a smooth variety and $\ua$, $\ub\subseteq\cO_X$ are nonzero sheaves
of ideals, then for every $\gamma\in\QQ_+$ we have
\begin{equation}
\cI(X, \gamma\cdot (\ua+\ub))\subseteq \sum_{\alpha+\beta=\gamma}
\cI(X, \alpha\cdot\ua)\cdot\cI(X, \beta\cdot\ub).
\end{equation}
\end{theorem} 

Note that the sum has finitely many distinct terms.
The above statement admits a generalization to the case of two graded systems
of ideals. Recall that a graded system of ideals $\ua_{\bullet}
=(\ua_m)_{m\geq 0}$ on $X$
is a sequence of nonzero ideals such that $\ua_0=\cO_X$ and 
$\ua_p\cdot\ua_q\subseteq\ua_{p+q}$, for every $p$ and $q$. It is proved in
\cite{els} that the set $\{\cI(X, \alpha/q\cdot\ua_{pq})\}_{q\geq 1}$
has a unique maximal element, the asymptotic multiplier ideal 
$\cI(X, \alpha\cdot\|\ua_p\|)$.

Suppose now that we have two graded systems of ideals $\ua_{\bullet}$
and $\ub_{\bullet}$. Their sum $\uc_{\bullet}=\ua_{\bullet}
+\ub_{\bullet}$, defined by $\uc_m=\sum_{i+j=m}\ua_i\ub_j$, 
is again a graded system of ideals. With these definitions we have the
 following

\begin{theorem}\label{asymptotic}
Let $X$ be a smooth variety and $\ua_{\bullet}$ and $\ub_{\bullet}$
two graded systems of ideals on $X$ and $\uc_{\bullet}$ their sum.
 For every $\gamma\in\QQ_+$ and every
$p\geq 1$, we have
\begin{equation}
\cI(X,\gamma\cdot \|\uc_p\|)\subseteq \sum_{\alpha+\beta=\gamma}
\cI(X,\alpha\cdot \|\ua_p\|)\cdot\cI(X,\beta\cdot \|\ub_p\|).
\end{equation}
\end{theorem}

Note that Theorem~\ref{one_ideal} can be obtained from Theorem~\ref{asymptotic}
by taking the systems $\ua_{\bullet}$ and $\ub_{\bullet}$ to be given
by the powers of $\ua$ and $\ub$, respectively.

\smallskip

As an application of Theorem~\ref{one_ideal}, in the second part of the paper
 we show that general
 multiplier ideals can be approximated at each point 
by multiplier ideals associated
to zero dimensional ideals (see Theorem~\ref{approximation} for the
precise statement). This is then used to compare the multiplier
ideals associated to a scheme in different embeddings. For example,
we show that if $\ua\subseteq\cO_X$ and $\ub\subseteq\cO_Y$ are such that
$V(\ua)\simeq V(\ub)$, then this isomorphism maps 
$\cI(X, \alpha\cdot\ua)\cdot\cO_{V(\ua)}$
to $\cI(Y, \beta\cdot\ub)\cdot\cO_{V(\ub)}$ if $\dim\,X-\alpha=\dim\,Y-\beta$.

\bigskip

We give below the idea of the proof of the main results. For simplicity,
we consider only the case of Theorem~\ref{one_ideal}. 
 The proof of
Theorem~\ref{asymptotic} follows from a similar, but more technical statement
which can be proved in an analogous way (see Theorem~\ref{main}).

The first step is to use the Restriction theorem for the canonical embedding
$\Delta\,:\,X\hookrightarrow X\times X$ to reduce the statement of
 Theorem~\ref{one_ideal} to a result on $X\times X$. This is the particular
case $X=Y$ in the following

\begin{theorem}\label{equality}
Let $X$ and $Y$ be smooth varieties and $\ua\subseteq \cO_X$
and $\ub\subseteq \cO_Y$ nonzero sheaves of ideals. If 
$p\,:\,X\times Y\longrightarrow X$ and $q\,:\,X\times Y\longrightarrow
Y$ are the canonical projections, then for every $\gamma\in\QQ_+$
we have
\begin{equation}\label{product}
\cI\left(X\times Y, \gamma\cdot (p^{-1}\ua+q^{-1}\ub)\right)
=\sum_{\alpha+\beta
=\gamma}p^{-1}\cI(X,\alpha\cdot\ua)\cdot q^{-1}\cI(Y,\beta\cdot\ub).
\end{equation}
\end{theorem}

The next step is to reduce Theorem~\ref{equality}
by taking log resolutions to the case
when $\ua$ and $\ub$ are ideals defining divisors with simple normal crossings.
 However,
as it stands, the righthand side of equation~\ref{product} does not behave
well with respect to push-forward. Therefore we first prove a lemma
showing that in fact, with the above notation, we have
$$\sum_{\alpha+\beta=\gamma}p^{-1}\cI(X,\alpha\cdot\ua)\cdot q^{-1}
\cI(Y,\beta\cdot\ub)=\bigcap_{\alpha+\beta=\gamma}\left
(p^{-1}\cI(X,\alpha\cdot\ua)
+q^{-1}\cI(Y,\beta\cdot\ub)\right).$$

Using this expression, we can 
reduce ourselves to the case of divisors with simple
normal crossings. Note, however that $p^{-1}\ua+q^{-1}\ub$ has codimension
two. On the other hand, 
locally in the \'{e}tale topology $\ua$ and $\ub$ are monomial ideals and
therefore so is $p^{-1}\ua+q^{-1}\ub$. The equality in Theorem~\ref{equality}
follows now using the explicit description of multiplier ideals of
monomial ideals due to Howald (see \cite{howald}).

\smallskip

\subsection{Acknowledgements}
We are grateful to Lawrence Ein and Rob Lazarsfeld for their 
encouragement and for their comments on this work. In particular, the 
extension of Theorem~\ref{one_ideal} to asymptotic multiplier ideals was
suggested to us by Rob Lazarsfeld.

\section{The multiplier ideals of a sum of ideals}

 We work over the field of complex numbers. All sheaves of ideals
are assumed to be 
quasi-coherent. The basic results on multiplier ideals that will be
used can be found in \cite{ein}, \cite{del} and \cite{els}
(see also \cite{lazarsfeld} for a thorough presentation of the subject).

\smallskip

The following theorem is the main technical
 result of this section. It easily implies
Theorems~\ref{one_ideal} and~\ref{asymptotic}.

\begin{theorem}\label{main}
Let $X$ and $Y$ be smooth varieties and $p\,:\,X\times Y\longrightarrow X$
and $q\,:\,X\times Y\longrightarrow Y$ the canonical projections.
Suppose that $\ua_i\subseteq\cO_X$, $\ub_i\subseteq\cO_Y$ are nonzero sheaves
of ideals for $1\leq i\leq m$ or $i=n$. If for every $1\leq i\leq m$
we have $i\mid n$ and $\ua_i^{n/i}\subseteq \ua_n$ and $\ub_i^{n/i}
\subseteq \ub_n$, then
\begin{equation}
\cI\left(X\times Y,
{\gamma}/m\cdot\left(p^{-1}\ua_m+\sum_{i=1}^{m-1}p^{-1}\ua_i\cdot 
q^{-1}\ub_{m-i}+q^{-1}\ub_m\right)\right)\subseteq
\end{equation}
$$\sum_{\alpha+\beta=\gamma}
p^{-1}\cI(X, {\alpha}/n\cdot\ua_n)
\cdot q^{-1}\cI(Y, {\beta}/n
\cdot\ub_n),$$
for every $\gamma\in\QQ_+$.
\end{theorem}

We prove first the following lemma.

\begin{lemma}\label{equiv}
Let $X$ and $Y$ be smooth varieties, $p\,:\,X\times Y\longrightarrow X$
and $q\,:\,X\times Y\longrightarrow Y$ the canonical projections,
and $\ua\subseteq\cO_X$ and $\ub\subseteq\cO_Y$ nonzero sheaves of ideals.
For every rational number $\gamma\geq 0$, we have
$$\sum_{\alpha+\beta=\gamma}p^{-1}\cI(X,\alpha\cdot\ua)\cdot q^{-1}
\cI(Y,\beta\cdot\ub)=\bigcap_{\alpha+\beta=\gamma}\left
(p^{-1}\cI(X,\alpha\cdot\ua)
+q^{-1}\cI(Y,\beta\cdot\ub)\right).$$
\end{lemma}

\begin{proof}
In order to prove ``$\subseteq$'', we have to show that if
$\alpha+\beta=\alpha_1+\beta_1=\gamma$, then
$$p^{-1}\cI(X,\alpha\cdot\ua)\cdot q^{-1}\cI(Y,\beta\cdot\ub)
\subseteq p^{-1}\cI(X,\alpha_1\cdot\ua)+q^{-1}\cI(Y,\beta_1\cdot\ub).$$
It is clear that we must have either $\alpha\geq\alpha_1$ or
 $\beta\geq \beta_1$. In the first case we have $\cI(X,\alpha\cdot\ua)
\subseteq\cI(X,\alpha_1\cdot\ua)$ and the above inclusion follows. The other
case is similar.

In order to prove the reverse inclusion, we may assume that $X$ and $Y$
are affine and let $R=\cO(X)$  and $S=\cO(Y)$. We  identify the
multiplier ideals with their global sections.

We choose bases for $R$ and $S$ which are compatible with the
filtrations induced by 
the respective multiplier ideals,
as follows. Note that the set of multiplier ideals of
$\ua$ with coefficient $p$, $0\leq p\leq \gamma$ is finite and
$$\cI(X, p\cdot\ua)\subseteq\cI(X, p'\cdot\ua),$$
if $p>p'$. Therefore we can choose index sets $I_{\lambda}$, possibly empty,
for $0\leq\lambda\leq\gamma$, and elements $e_{\lambda\mu}\in R$,
for $0\leq\lambda\leq\gamma$ and $\mu\in I_{\lambda}$ such that for
every $p\leq\gamma$, a basis over $\CC$ for
$\cI(X,p\cdot\ua)$ is given by
$\{e_{\lambda\mu}\,\mid\,\lambda\leq\gamma-p,\, \mu\in I_{\lambda}\}$.

We consider an analogous set of elements $f_{\lambda\mu}\in S$,
with $0\leq\lambda\leq\gamma$ and $\mu\in J_{\lambda}$,
satisfying the corresponding property with respect to the multiplier ideals of
$\ub$.

A basis in $p^{-1}\cI(X,\alpha\cdot\ua)\cdot q^{-1}\cI(Y,\beta\cdot b)$
is given by 
$$\{e_{\lambda\mu}\otimes f_{\lambda'\mu'}\,\mid\,
\lambda\leq\gamma-\alpha,\,\lambda'\leq\gamma-\beta,\,\mu\in I_{\lambda},\,
\mu'\in J_{\lambda'}\}.$$

Therefore a basis in $\sum_{\alpha+\beta=\gamma}p^{-1}\cI(X,\alpha\cdot\ua)
\cdot q^{-1}\cI(Y,\beta\cdot\ub)$ is given by
$$\{e_{\lambda\mu}\otimes f_{\lambda'\mu'}\,\mid\,\lambda+\lambda'
\leq\gamma,\, \mu\in I_{\lambda},\, \mu'\in J_{\lambda'}\}.$$

It is enough to prove that if 
$$e_{\lambda\mu}\otimes f_{\lambda'\mu'}\in\bigcap_{\alpha+\beta=\gamma}
\left(p^{-1}\cI(X,\alpha\cdot\ua)+ q^{-1}\cI(Y,\beta\cdot\ub)\right),$$
then $\lambda+\lambda'\leq\gamma$. Indeed, the above intersection
has a basis given by a subset of $\{e_{\lambda\mu}\otimes
e_{\lambda'\mu'}\}_{\lambda,\mu,\lambda',\mu'}$, because so has
each member of the intersection.

For every $\alpha$ and $\beta$ such that $\alpha+\beta=\gamma$, we must have
either $\lambda\leq\gamma-\alpha$ or $\lambda'\leq\gamma-\beta$.
Therefore for every $0\leq\alpha\leq\gamma$, we have either $\lambda\leq 
\gamma-\alpha$ or $\lambda'\leq\alpha$. This gives
$\lambda'\leq\gamma-\lambda$ and finishes the proof of the lemma.
\end{proof}

The proof of Theorem~\ref{main} will be based on reduction to the
case of monomial ideals. Therefore we first treat this special case
in the following lemma.

\begin{lemma}\label{monomial}
The statement of Theorem~\ref{main} is true if $X=\AA^r$, $Y=\AA^s$
and $\ua_i\subseteq \CC[X]=\CC[X_1,\ldots,X_r]$ and $\ub_i\subseteq \CC[Y]
=\CC[Y_1,\ldots,Y_s]$ are monomial ideals for all $i$.
\end{lemma}

\begin{proof}
By Lemma~\ref{equiv}, it is enough to prove that
\begin{equation}\label{eqn}
\cI\left(\AA^r\times \AA^s,\gamma/m\cdot\left(p^{-1}\ua_m+
\sum_{i=1}^{m-1}p^{-1}\ua_i\cdot q^{-1}
\ub_{m-i}+q^{-1}\ub_m\right)\right)\subseteq
\end{equation}
$$\bigcap_{\alpha+\beta=\gamma}
\left(p^{-1}\cI(\AA^r,\alpha/n\cdot\ua_n)+
 q^{-1}\cI(\AA^s,\beta/n\cdot\ub_n)\right).$$

If there is $i$ such that $\ua_i=\CC[X]$, then $\ua_n=\CC[X]$ and therefore
$\cI(X, \alpha\cdot\ua_n)=\CC[X]$, for every $\alpha$. It follows that the
right hand side of equation~(\ref{eqn})
 is $\CC[X,Y]$ and the inclusion is obvious.

We may therefore assume that $\ua_i\neq \CC[X]$ for every $i$, and by symmetry,
that $\ub_i\neq \CC[Y]$, for every $i$. Suppose that for some $u\in\NN^r$
and $v\in\NN^s$, we have
$$X^uY^v\in\cI\left(\AA^r\times \AA^s,
\gamma/m\cdot \left(p^{-1}\ua_m+\sum_{i=1}^{m-1}
p^{-1}\ua_i\cdot q^{-1}\ub_{m-i}+q^{-1}\ub_m\right)\right),$$
but for some $\alpha$, $\beta\in\QQ_+$ with $\alpha+\beta=\gamma$,
we have
 $X^u\not\in\cI(\AA^r,\alpha/n\cdot\ua_n)$ and
 $Y^v\not\in\cI(\AA^s,\beta/n\cdot\ub_n)$.

We use Howald's description for multiplier ideals
of monomial ideals 
 in \cite{howald}. It says that if $I\subsetneqq \CC[X]$
is a nonzero monomial ideal and $P_I\subset \RR^r$ is the convex hull of
$\{w\in\NN^r\,\mid\,X^w\in I\}$, then for every $c>0$,
$$\cI({\bf A}^r,c\cdot I)=\left(X^w\,\mid\,w+e\in {\rm Int}\,(c\cdot P_I)
\right).$$
Here $e$ denotes the unit vector $(1,\ldots,1)\in\NN^r$.

Our hypothesis on $u$ and $v$ implies that $u+e\not\in {\rm Int}\,(\alpha/n
\cdot P_{\ua_n})$ and $v+f\not\in {\rm Int}\,(\beta/n\cdot P_{\ub_n})$
(where $f=(1,\ldots,1)\in\NN^s$). This means that there are linear maps
$\phi\,:\,\RR^r\longrightarrow\RR$
 and $\psi\,:\,\RR^s\longrightarrow\RR$ such that
$\phi(w_1)\geq 1$ if $X^{w_1}\in\ua_n$ and $\psi(w_2)\geq 1$ if
$Y^{w_2}\in\ub_n$, but $\phi(u+e)\leq \alpha/n$ and
$\psi(v+f)\leq \beta/n$. Therefore we have
$\phi(u+e)+\psi(v+f)\leq\gamma/n$.

If $X^{w_1}Y^{w_2}\in p^{-1}\ua_i\cdot q^{-1}\ub_{m-i}$,
since $\ua_i^{n/i}\subseteq\ua_n$ and $\ub_i^{n/i}\subseteq \ub_n$, 
we get
$(n/i)w_1\in P_{\ua_n}$ and $(n/(m-i))w_2\in P_{\ub_n}$. We deduce
that $\phi(w_1)+\psi(w_2)\geq i/n+(m-i)/n=m/n$. If $X^{w_1}Y^{w_2}\in
p^{-1}\ua_m$, then $\phi(w_1)\geq m/n$ and $\psi(w_2)\geq 0$ and we have
 analogous inequalities if $X^{w_1}Y^{w_2}\in q^{-1}\ub_m$.

This shows that the linear 
map $\rho\,:\,\RR^r\times\RR^s\longrightarrow
\RR$, given by $\rho(w_1,w_2)=\phi(w_1)+\psi(w_2)$ has the property that
$(n/m)\rho(w_1,w_2)\geq 1$ if
 $$X^{w_1}Y^{w_2}\in
p^{-1}\ua_m+\sum_{i=1}^{m-1}p^{-1}\ua_i\cdot q^{-1}\ub_{m-i}+q^{-1}\ub_m.$$
Since we have
 $$X^uY^v\in\cI\left(\AA^r\times \AA^s,\gamma/m\cdot
\left(p^{-1}\ua_m+\sum_{i=1}^{m-1}p^{-1}\ua_i
\cdot q^{-1}\ub_{m-i}+q^{-1}\ub_m\right)\right),$$
one more application of Howald's theorem gives
$n/m\left(\phi(u+e)+\psi(v+f)\right)>\gamma/m$, a contradiction.
\end{proof}

We can give now the proof of the general case.

\begin{proof}[Proof of Theorem~\ref{main}]
Let $f\,:\,X'\longrightarrow X$ and
$g\,:\,Y'\longrightarrow Y$ be log resolutions
for all pairs $(X,\ua_i)$ and $(Y, \ub_i)$, and also
for $\left(X,\sum_i\ua_i\right)$ and $\left(Y,\sum_i\ub_i\right)$,
 respectively.

Let $\ua'_i=f^{-1}\ua_i$ and $\ub'_i=g^{-1}\ub_i$.
If $p'\,:\,X'\times Y'\longrightarrow X'$ and
$q'\,:\,X'\times Y'\longrightarrow Y'$ are the canonical projections,
then
we use the notation
$$\ur=p^{-1}\ua_m+\sum_{i=1}^{m-1}p^{-1}\ua_i
\cdot q^{-1}\ub_{m-i}+q^{-1}\ub_m,$$
$$\ur'={p'}^{-1}\ua'_m+\sum_{i=1}^{m-1}{p'}^{-1}\ua'_i
\cdot {q'}^{-1}\ub'_{m-i}+{q'}^{-1}\ub'_m.$$
If $h=(f,g)\,:\,X'\times Y'\longrightarrow X\times Y$, then
$\ur'=h^{-1}\ur$.

Suppose first that the assertion of the theorem
is true for $X'$, $Y'$, $\{\ua'_i\}_i$ and $\{\ub'_i\}_i$.
The change of variable formula for multiplier ideals gives
$$\cI(X\times Y, \gamma/m\cdot\ur)
=h_*\left(\cI(X'\times Y',\gamma/m\cdot\ur')\otimes
\cO(K_{X'\times Y'/X\times Y})\right).$$
Using also Lemma~\ref{equiv}, we deduce
$\cI(X\times Y,\gamma/m\cdot\ur)\subseteq$
$$h_*\left(\bigcap_{\alpha+\beta=\gamma}
\left({p'}^{-1}\cI(X',\alpha/n\cdot\ua'_n)+
{q'}^{-1}\cI(Y',\beta/n\cdot\ub'_n)\right)\otimes\cO(K_{X'\times Y'
/X\times Y})\right)$$
$$=\bigcap_{\alpha+\beta=\gamma}h_*\left(
\left({p'}^{-1}\cI(X', \alpha/n\cdot\ua'_n)+
{q'}^{-1}\cI(Y', \beta/n\cdot\ub'_n)\right)\otimes\cO(K_{X'\times Y'/
X\times Y})\right).$$
Indeed, the sets $\{\cI(X',\alpha\cdot\ua'_n)\mid\alpha\leq\gamma/n\}$
and $\{\cI(Y',\beta\cdot\ua'_n)\mid\beta\leq\gamma/n\}$ are finite,
so that the above intersection has finitely many distinct terms,
and therefore commutes with push-forward.

Note that if $\cF'$, $\cF''\subseteq \cF$ are quasicoherent subsheaves
of $\cO_{X'\times Y'}$--modules and if $R^1h_*(
\cF'\cap \cF'')=0$, then we have $h_*(\cF'+\cF'')=h_*\cF'+h_*\cF''$.
Indeed, this follows by applying $h_*$ to the
exact sequence
$$0\longrightarrow\cF'\cap\cF''\longrightarrow
\cF'\oplus\cF''\longrightarrow\cF'+\cF''\longrightarrow 0.$$

Note that
$$R^1h_*\left(\left({p'}^{-1}\cI(X',\alpha/n\cdot\ua'_n)\cap {q'}^{-1}
\cI(Y',\beta/n\cdot\ub'_n)\right)\otimes\cO
(K_{X'\times Y'/
X\times Y})\right)=$$
$$R^1h_*\left({p'}^{-1}\left(\cI(X',\alpha/n\cdot\ua'_n)\otimes
\cO(K_{X'/X})\right)\otimes {q'}^{-1}\left(\cI(Y',\beta/n\cdot\ub'_n)
\otimes
\cO(K_{Y'/Y})\right)\right)$$
vanishes.
This follows by applying the K\"{u}nneth formula and the
Local Vanishing theorem (see \cite{ein} 1.4)
which gives $R^1f_*(\cI(X',\alpha/n\cdot\ua'_n)\otimes
\cO(K_{X'/X}))=0$ and $R^1g_*(\cI(Y',\beta/n\cdot\ub'_n)
\otimes\cO(K_{Y'/Y}))=0$.

Using the fact that $f_*(\cO(K_{X'/X}))=\cO_X$
and $g_*(\cO(K_{Y'/Y}))=\cO_Y$, one more application of
the K\"{u}nneth formula and of the change of variable formula
for the multiplier ideals gives:
$$h_*\left({p'}^{-1}\cI(X',\alpha/n\cdot\ua'_n)
\otimes\cO(K_{X'\times Y'/X\times Y})\right)
=p^{-1}\cI(X,\alpha/n\cdot\ua_n),$$
$$h_*\left({q'}^{-1}\cI(Y',\beta/n\cdot\ub'_n)\otimes
\cO(K_{X'\times Y'/X\times Y})\right)=q^{-1}\cI(Y,\beta/n\cdot\ub_n).$$

Putting everything together, we get via Lemma~\ref{equiv}
the statement of the theorem. To finish the proof, it is therefore
enough to consider the case when all $\ua_i$ and $\ub_i$
are ideals defining effective divisors on $X$ and $Y$, respectively,
whose union has simple normal crossings.

Since the statement of the theorem is local in $X$ and $Y$, we may
assume that we have \'{e}tale morphisms $\phi\,:\,X\longrightarrow
\AA^r$ and $\psi\,:\,Y\longrightarrow\AA^s$ whose images contain
the origins in the respective affine spaces, and principal monomial ideals
$\widetilde{{\ua}_i}$ and $\widetilde{{\ub}_i}$
 such that $\ua_i=\phi^{-1}\widetilde{{\ua}_i}$
and $\ub_i=\psi^{-1}\widetilde{{\ub}_i}$, for all $i$.

Since $\phi$ and $\psi$ are \'{e}tale, the hypothesis implies
$\widetilde{{\ua}_i}^{n/i}\subseteq\widetilde{\ua}_n$ and 
$\widetilde{{\ub}_i}^{n/i}\subseteq\widetilde{{\ub}_n}$,
 for all $1\leq i\leq n$.
Moreover, taking multiplier ideals commutes with the pull-back by
\'{e}tale morphisms, so that we
can reduce the theorem to the case of monomial ideals, when it follows
from Lemma~\ref{monomial}. 
\end{proof}

As in \cite{del}, 
we can use the Restriction theorem to deduce from
 Theorem~\ref{main}
a property of families of ideals on the same variety.

\begin{corollary}\label{same_var}
Let $X$ be a smooth variety and $\ua_i\subseteq\cO_X$, $\ub_i\subseteq\cO_X$
nonzero sheaves of ideals, with $1\leq i\leq m$ or $i=n$. If for every
$1\leq i\leq m$, we have $i\mid n$ and $\ua_i^{n/i}\subseteq\ua_n$,
$\ub_i^{n/i}\subseteq\ub_n$, then for every 
$\gamma\in\QQ_+$, we get
$$\cI(X,\gamma/m\cdot(\ua_m+\sum_{i=1}^{m-1}
\ua_i\cdot\ub_{m-i}+\ub_m))
\subseteq\sum_{\alpha+\beta=\gamma}\cI(X,\alpha/n\cdot\ua_n)
\cdot\cI(X,\beta/n\cdot\ub_n).$$
\end{corollary}

\begin{proof}
Consider the diagonal embedding $X\hookrightarrow X\times X$. 
If $p\,:\,X\times X\longrightarrow X$ and $q\,:\,
X\times X\hookrightarrow X$ are the projections on the first and,
respectively, the second component, let
$$\ur=p^{-1}\ua_m+\sum_{i=1}^mp^{-1}\ua_i\cdot q^{-1}\ub_{m-i}
+q^{-1}\ub_m.$$
Note that we have $\ur\cdot\cO_X=\ua_m+\sum_{i=1}^{m-1}\ua_i\cdot\ub_{m-i}
+\ub_m$.

We clearly have $X\not\subseteq\,V(\ur)$, so that by the
Restriction theorem (see \cite{ein} 2.1)
we deduce
$$\cI\left(X,\gamma/m\cdot\left(\ua_m+\sum_{i=1}^{m-1}
\ua_i\cdot\ub_{m-i}+\ub_m\right)\right)
\subseteq\cI(X\times X,\gamma/m\cdot\ur)\cdot\cO_X.$$

On the other hand, Theorem~\ref{main} gives
$$\cI(X\times X,\gamma/m\cdot\ur)\cdot\cO_X
\subseteq \left(\sum_{\alpha+\beta=\gamma}p^{-1}\cI(X,\alpha/n\cdot\ua_n)
\cdot q^{-1}\cI(X,\beta/n\cdot\ub_n)\right)\cdot\cO_X$$
$$=\sum_{\alpha+\beta=\gamma}\cI(X,\alpha/n\cdot\ua_n)\cdot
\cI(X,\beta/n\cdot\ub_n).$$
The above inclusions imply the statement of the corollary.
\end{proof}

We can give now the proofs of the statements announced in the Introduction.

\begin{proof}[Proof of Theorem~\ref{asymptotic}]
Using the fact that 
$$\cI(X, \alpha\cdot\|\ua_p\|)=\cI(X,p\alpha\cdot \|\ua_1\|)$$
and similar equalities for $\ub_{\bullet}$ and $\uc_{\bullet}$,
we  reduce immediately to the case
$p=1$. By definition, we have to prove that for every $m\geq 1$, we have
$$\cI(X, \gamma/m\cdot\uc_m)\subseteq\sum_{\alpha+\beta=\gamma}
\cI(X, \alpha\cdot \|\ua_1\|)\cdot\cI(X, \beta\cdot \|\ub_1\|).$$

If $n$ is a positive integer such that for every $1\leq i\leq m$ we have
$i\mid n$, then we can apply
Corollary~\ref{same_var} to get
$$\cI(X, \gamma/m\cdot\uc_m)\subseteq
\sum_{\alpha+\beta=\gamma}\cI(X, \alpha/n\cdot\ua_n)\cdot\cI(X,
\beta/b\cdot\ub_n).$$

On the other hand we have by definition $\cI(X, \alpha/n\cdot\ua_n)
\subseteq\cI(X,\alpha\cdot \|\ua_1\|)$ and a similar inclusion for 
$\ub_{\bullet}$. This proves the statement of the corollary.
\end{proof}

\smallskip

\begin{proof}[Proof of Theorem~\ref{one_ideal}]
This is precisely the statement of Corollary~\ref{same_var} in the case
$m=n=1$.
\end{proof}

We give now the proof of Theorem~\ref{equality}. Recall that it says that
in the particular case $m=n=1$, the inclusion in Theorem~\ref{main}
becomes equality.

\begin{proof}[Proof of Theorem~\ref{equality}]
The proof of Theorem~\ref{main} applies word by word in this case
if we know that we have equality for monomial ideals.
Therefore we may assume that $X=\AA^r$, $Y=\AA^s$
and that $\ua\subseteq \CC[X]=\CC[X_1,\ldots,X_r]$ and
$\ub\subseteq \CC[Y]=\CC[Y_1,\ldots,Y_s]$ are monomial ideals.
We have to prove that if $\alpha+\beta=\gamma$, then
$$p^{-1}\cI(\AA^r,\alpha\cdot\ua)\cdot q^{-1}\cI(\AA^s, \beta\cdot\ub)
\subseteq\cI\left(\AA^r\times \AA^s,
 \gamma\cdot(p^{-1}\ua+ q^{-1}\ub)\right).$$

If $\ua=\CC[X]$ or $\ub=\CC[Y]$, then the right hand side of the above
inclusion is $\CC[X,Y]$, and the statement is trivial. Suppose
 therefore that we are in none of these cases. Moreover,
if $\alpha=0$, then it is easy to see from the definition of
 multiplier ideals that $q^{-1}\cI(\AA^s,\gamma\cdot\ub)
=\cI(\AA^r\times \AA^s, \gamma\cdot q^{-1}\ub)$ (note that Proposition~2.2
 in \cite{del} and its extension to the case of ideals
 give a more general 
statement).
 In this case we get the above inclusion since
$\cI(\AA^r\times \AA^s, \gamma\cdot q^{-1}\ub)\subseteq
\cI\left(\AA^r\times \AA^s, \gamma\cdot (p^{-1}\ua+q^{-1}\ub)\right)$.
Therefore we may assume that $\alpha>0$, and by symmetry,
also that $\beta>0$.

We use again the description in \cite{howald} for multiplier ideals
of monomial ideals. First, this description shows that these ideals
are generated by monomials.

Suppose that we have $X^uY^v\in p^{-1}\cI(\AA^r, \alpha\cdot\ua)\cdot
q^{-1}\cI(\AA^s, \beta\cdot\ub)$. 
Let $\phi\,:\,\RR^r\times\RR^s\longrightarrow
\RR$ be a linear map such that $\phi(w_1,0)\geq 1$ if
$X^{w_1}\in\ua$ and $\phi(0,w_2)\geq 1$ if
$Y^{w_2}\in\ub$. By \cite{howald}, it is enough 
to prove that for every such $\phi$
we have $\phi(u,v)>\gamma$.

On the other hand, since $X^u\in\cI(\AA^r, \alpha\cdot\ua)$ we get
$\phi(u,0)>\alpha$. Similarly, since $Y^v\in\cI(\AA^s, \beta\cdot\ub)$
we get $\phi(0,v)>\beta$. This implies that $\phi(u,v)>\gamma$.
\end{proof}

\begin{remark}
If we make the convention that $\cI(X, \gamma\cdot (0))$
is equal to $\cO_X$ if $\gamma=0$ and $(0)$ otherwise, then
the formula in Theorem~\ref{equality} is still valid
if $\ua=(0)$ or $\ub=(0)$. Indeed, if for example 
$\ua=(0)$, then the formula in Theorem~\ref{equality} says that
$$\cI(X\times Y, \gamma\cdot q^{-1}\ub)=q^{-1}\cI(Y, \gamma\cdot
\ub).$$
But as we have mentioned, this is a particular case of the results in 
\cite{del}.
\end{remark}

\section{Invariance of multiplier ideals}

We start by showing that Theorem~\ref{one_ideal} can be used to approximate
arbitrary multiplier ideals by multiplier ideals corresponding to 
zero dimensional ideals.

Fix a smooth variety $X$ with $\dim\,X=n$
and $\ua\subseteq\cO_X$ a nonzero sheaf of ideals
on $X$.
For a (closed) point $x\in X$, let $\um_x\subseteq\cO_x$ be the
ideal defining that point. If $l\leq 0$ is an integer, we put 
$\um_x^l=\cO_X$.
Recall that $[\cdot]$ denotes the integral part function.

\begin{proposition}\label{approximation}
With the above notation, for every $x\in X$, every integer $p\geq 1$,
and every $\gamma\in\QQ_+$, $\epsilon\in\QQ_+^*$, 
 we have
\begin{equation}
\cI\left(X, (\gamma+\epsilon)\cdot (\ua+\um_x^p)\right)\subseteq
\cI(X, \gamma\cdot\ua)+\um_x^{[p\epsilon]-n+1}.
\end{equation}
\end{proposition}

\begin{proof}
We apply Theorem~\ref{one_ideal} to the sheaves of ideals $\ua$
and $\um_x^p$ and to $\gamma+\epsilon$. Note that the multiplier
ideals of $\um_x^p$ are given by
$$\cI(X, \beta\cdot\um_x^p)=\um_x^{[p\beta]-n+1},$$
for every $p\geq 1$.
We therefore obtain
$$\cI\left(X, (\gamma+\epsilon)\cdot(\ua+\um_x^p)\right)\subseteq
\sum_{\alpha+\beta=\gamma+\epsilon}\cI(X, \alpha\cdot\ua)
\cdot \cI(X, \beta\cdot\um_x^p)\subseteq$$
$$\subseteq\cI(X, \gamma\cdot\ua)+\sum_{\beta>\epsilon}
\cI(X, \beta\cdot\um_x^p)=\cI(X, \gamma\cdot\ua)
+\cI(X, \epsilon\cdot\um_x^p).$$

The statement of the proposition now follows from this
and the formula for the multiplier ideals of $\um_x^p$.
\end{proof}

\begin{remark}
With the notation in the above proposition, let $I=
\cI(X, \gamma\cdot\ua)$. By the semicontinuity property of
multiplier ideals with respect to the parameter, we can find
 $\epsilon\in\QQ_+^*$ such that
 $\cI\left(X, (\gamma+\epsilon)\cdot\ua\right)=I$.
We deduce from this and Proposition~\ref{approximation} that

$$I+\um_x^{[p\gamma+p\epsilon]-n+1}\subseteq 
\cI\left(X,(\gamma+\epsilon)\cdot (\ua+\um_x^p)\right)
\subseteq I+\um_x^{[p\epsilon]-n+1}.$$
Since multiplier ideals are integrally closed, we deduce that
the integral closure of $I+\um_x^{[p\gamma+p\epsilon]-n+1}$
is contained in $I+\um_x^{[p\epsilon]-n+1}$. 

 As pointed to us by L.Ein and R. Lazarsfeld,
this can be considered as an effective version of a theorem
of Delfino and Swanson 
(see \cite{ds}) for the case when $I$ is a multiplier ideal (note however that
the result in \cite{ds} holds in arbitrary excellent rings).
\end{remark}

\bigskip

We first apply this to study the relation between multiplier ideals
on $X$ and $Y$, when $X$ is a subvariety of $Y$. We use the convention
that if $\ua\subseteq\cO_X$ is a nonzero sheaf of ideals on a smooth variety
 $X$, then $\cI(X, \gamma\cdot\ua)=\cO_X$, for every $\gamma<0$.

\begin{proposition}\label{subvariety}
Let $Y$ be a smooth variety and $X\subset Y$ a closed smooth subvariety,
with ${\rm codim}\,(X/Y)=r$. If $\ua_X\subset\cO_Y$ is the sheaf of ideals
defining $X$ and $\ub$ is a sheaf of ideals such that $\ua_X\subsetneqq\ub$,
then 
$$\cI(X,\gamma\cdot\ub/\ua_X)=\cI\left(Y, (\gamma+r)\cdot\ub\right)
\cdot\cO_X,$$
for every $\gamma\in\QQ$.
\end{proposition}

\begin{proof}
We consider first the case when
 $$X=\AA^n\hookrightarrow Y=X\times
\AA^r,$$
is defined by the vanishing of the last $r$ coordinates.
Let $p\,:\,X\times\AA^r\longrightarrow X$
and $q\,:\,X\times\AA^r\longrightarrow\AA^r$ be the canonical
projections.

 If $\um_0$ is the ideal defining the origin in $\AA^r$,
then $\ua_X=q^{-1}\um_0$. Moreover, we have
$\ub= p^{-1}(\ub/\ua_X)+q^{-1}\um_0$. We know that 
$\cI(\AA^r, \beta\cdot\um_0)$ is equal to $\um_0^{[\beta]-r+1}$, if
$\beta\geq r-1$, and it is equal with $\cO_{\AA^r}$, otherwise.

We may assume that $\gamma+r\geq 0$ because otherwise the statement
of the proposition is trivial.
Theorem~\ref{equality} gives
$$\cI\left(Y, (\gamma+r)\cdot\ub\right)=
\sum_{\alpha+\beta=\gamma+r}p^{-1}\cI(X, \alpha\cdot \ub/\ua_X)
\cdot q^{-1}\cI(\AA^r, \beta\cdot\um_0).$$

Since we have $q^{-1}\cI(\AA^r, \beta\cdot\um_0)\cdot\cO_X=(0)$
if $\beta\geq r$, we deduce that
$$\cI\left(Y, (\gamma+r)\cdot\ub\right)\cdot\cO_X
=\sum_{\gamma<\alpha\leq\gamma+r}\cI(X,\alpha\cdot\ub/\ua_X)
=\cI(X, \gamma\cdot\ub/\ua_X),$$
which finishes the proof of this case.

\smallskip

We show now that if $\ub$ is a zero dimensional ideal, then we can reduce
the statement to the above case. Since the statement is local, we may
assume that ${\rm Supp}\,(\cO_Y/\ub)=\{x\}$, for some point $x\in X$ and
it is enough to check the equality in the proposition in an open
neighbourhood of $x$. Therefore we may assume that there is an 
\'{e}tale morphism $\phi\,:\,Y\longrightarrow\AA^{n+r}$
with $\phi(x)=0$ such that
$X=\phi^{-1}(\AA^n)$. Here we view $\AA^n$ embedded 
in $\AA^{n+r}$ as before.

Note that $\phi$ induces an isomorphism between the completions of
the local rings of $Y$ and $\AA^{n+r}$ at $x$ and $0$, respectively.
But $\dim\,(\cO_Y/\ub)_x=0$, so that there is an ideal
$\ub'\subset\cO_{\AA^{n+r}}$, such that $\ub=\phi^{-1}\ub'$.
Since construction of multiplier ideals commutes with pull-back by
\'{e}tale morphisms,
we deduce the proposition in the case of zero-dimensional ideals
from the case we have already proved.

\smallskip

To finish the proof of the proposition, we show how to deduce the
general case from that of zero dimensional ideals.
 Obviously it is enough to prove
that for every $x\in X$ we have equality after localizing at $x$:
$I\cdot\cO_{X,x}=J\cdot\cO_{X,x}$, where
$I:=\cI(X,\gamma\cdot\ub/\ua_X)$ and
$J:=\cI\left(Y, (\gamma+r)\cdot\ub\right)\cdot\cO_X$.

We fix $\epsilon\in\QQ_+^*$ such that $\cI(X, \gamma\cdot\ub/\ua_X)=
\cI\left(X, (\gamma+\epsilon)\cdot\ub/\ua_X\right)$
and $\cI\left(Y, (\gamma+r)\cdot\ub\right)
=\cI\left(Y, (\gamma+r+\epsilon)\cdot\ub\right)$.
 Let $\um_x$ be the ideal of the point
$x$ in $Y$.

 Using Proposition~\ref{approximation}
and the fact that we know the statement for
the zero-dimensional ideal
$\ub+\um_x^p$, for every $p\geq 1$, we get:
$$I\subseteq\cI\left(X, (\gamma+\epsilon)\cdot (\ub+\um_x^p)
/\ua_X\right)=
\cI\left(Y, (\gamma+\epsilon+r)\cdot (\ub+\um_x^p)\right)\cdot\cO_X$$
$$\subseteq
J+(\um_x/\ua_X)^{[p\epsilon]
-\dim\,Y+1}.$$
Since this is true for every $p\geq 1$, Krull's Intersection theorem
gives the inclusion $I\cdot\cO_{X,x}\subseteq J\cdot\cO_{X,x}$.

The reverse inclusion follows similarly:
$$J\subseteq
\cI\left(Y, (\gamma+\epsilon+r)\cdot (\ub+\um_x^p)\right)\cdot\cO_X$$
$$=\cI\left(X, (\gamma+\epsilon)\cdot(\ub+\um_x^p)/\ua_X\right)
\subseteq
J+(\um_x/\ua_X)^
{[p\epsilon]-\dim\,X+1},$$
for every $p\geq 1$, and we apply again Krull's Intersection theorem.
 This completes the proof of the proposition.
\end{proof}

We compare now multiplier ideals corresponding to a scheme in two different
arbitrary embeddings. More precisely, we show that 
the restrictions of the multiplier ideals
of $\ua$ to the subscheme defined by $\ua$ depend only on this subscheme
(they do not depend on the particular embedding into a smooth variety).

\begin{proposition}\label{invariance}
If $X_1$ and $X_2$ are smooth varieties and $\ua_1\subseteq\cO_{X_1}$,
$\ua_2\subseteq\cO_{X_2}$ are nonzero sheaves of ideals, then every isomorphism
$\phi\,:\,Y_1=V(\ua_1)\longrightarrow Y_2=V(\ua_2)$ induces isomorphisms
of ideals
$$\cI\left(X_1, (\gamma+\dim\,X_1)\cdot\ua_1\right)\cdot\cO_{Y_1}
\simeq\cI\left(X_2, (\gamma+\dim\,X_2)\cdot\ua_2\right)\cdot\cO_{Y_2},$$
for every $\gamma\in\QQ$.
\end{proposition}

\begin{proof}
If for example $\dim\,X_1=\dim\,X_2+s$, with $s\geq 1$, by replacing
$X_2$ with $X_2\times\AA^s$ and applying Proposition~\ref{subvariety},
we reduce ourselves to the case when $\dim\,X_1=\dim\,X_2=n$ and in this case 
it is enough to prove that for every $\gamma\in\QQ_+$, we get an
induced isomorphism:
\begin{equation}\label{inv}
\cI\left(X_1, \gamma\cdot\ua_1\right)\cdot\cO_{Y_1}
\simeq\cI\left(X_2, \gamma\cdot\ua_2\right)\cdot\cO_{Y_2}.
\end{equation}

As in the proof of Proposition~\ref{subvariety}, we first prove
the case $\dim\,Y_1=\dim\,Y_2=0$. To simplify the notation,
whenever there is no danger of confusion, we will
identify $Y_1$ with $Y_2$ via $\phi$ and denote it by $Y$. We may clearly
assume that the support of $Y$ consists of only one point $y\in Y$.

Let $r=\dim\,T_yY$. We pick a regular system of parameters
$x_1,\ldots, x_n$ for $X_1$ around $y$ such that $x_{r+1},\ldots,x_n$
are in the ideal of $Y_1$. After restricting to a suitable open
neighbourhood of $y$,
 this induces an \'{e}tale morphism
$\psi_1\,:\,X_1\longrightarrow\AA^n$ such that  
$Y_1\subseteq\psi_1^{-1}(\AA^r)$, where $\AA^r\subseteq\AA^n$
is defined by the vanishing of the last $n-r$ coordinates.
Moreover, since $\dim\,Y=0$, there is a subscheme $Y'_1\subseteq
\AA^r$ such that $Y_1=\psi_1^{-1}(Y'_1)$ and $\psi_1$ induces an isomorphism
$Y_1\simeq Y'_1$. 

We get a similar morphism $\psi_2\,:\,X_2\longrightarrow\AA^n$
with analogous properties. Using the fact that construction of multiplier
ideals commutes with pull-back by \'{e}tale morphisms we reduce the
equality in equation~(\ref{inv}) to the case when $X_1$ and $X_2$ are
both affine spaces. Moreover, using Proposition~\ref{subvariety}, we see
that we may assume that $r=n$.

In this case the isomorphism $\phi\,:Y_1\longrightarrow Y_2$ can
be lifted to a local ring homomorphism
$$\tilde{\phi}\,:\,{\cO}_{X_2,y}\longrightarrow
{\cO}_{X_1,y}.$$
Since $T_yY_i=T_yX_i$ for $i=1$, $2$, it follows that $\tilde{\phi}$
induces an isomorphism of the corresponding completion rings i.e.
it is \'{e}tale. By restricting further to neighbourhoods of $y$ in
$X_1$ and $X_2$, we may assume that $\tilde{\phi}$ is induced by an
\'{e}tale scheme morphism $X_1\longrightarrow X_2$. Using one more time
the invariance of multiplier ideals under pull-back for
\'{e}tale morphisms, we deduce equation~(\ref{inv})
 in the zero-dimensional case.

Suppose now that $Y$ has arbitrary dimension. It is enough to prove
that for every $y\in Y$, the analogue of equation~(\ref{inv})
holds for the images of those two ideals in $\cO_{Y,y}$. 
 Let us denote by 
$\um_1$ and $\um_2$ the ideals defining $y$ in $X_1$ and $X_2$
respectively. It is clear that we get induced isomophisms
$\phi_p\,:\,V(\ua_1+\um_1^p)\longrightarrow V(\ua_2
+\um_2^p)$, for every $p\geq 1$. 

We apply the statement in the case of the zero dimensional schemes
$Y_1^p=V(\ua_1+\um_1^p)$ and $Y_2^p=
V(\ua_2+\um_2^p)$. If $\epsilon\in\QQ_+^*$ is such that
$I(X_2, \gamma\cdot\ua_2)=\cI(X_2, (\gamma+\epsilon)\cdot\ua_2)$,
then applying also Proposition~\ref{approximation}, we get
$$\cI(X_1, \gamma\cdot\ua_1)\cdot\cO_{Y_1^p}\subseteq
\cI\left(X_1, (\gamma+\epsilon)\cdot (\ua_1+\um_1^p)\right)\cdot\cO_{Y_1^p}$$
$$\simeq\cI\left(X_2, (\gamma+\epsilon)\cdot (\ua_2+\um_2^p)\right)\cdot 
\cO_{Y_2^p}\subseteq \cI(X_2,\gamma\cdot\ua_2)\cdot\cO_{Y_2^p}
+\um_2^{[p\epsilon]-n+1}\cdot\cO_{Y_2^p}.$$

This shows that the image of $\cI(X_1, \gamma\cdot\ua_1)\cdot\cO_{Y_1,y}$
by the isomorphism induced by $\phi$ is contained in
$$\bigcap_{q\geq 1}\left(\cI(X_2, \gamma\cdot\ua_2)\cdot\cO_{Y_2,y}+
\um_2^q\cdot\cO_{Y_2,y}\right).$$
Krull's Intersection theorem gives now the inclusion $''\subseteq''$
in equation~(\ref{inv}). The reverse inclusion follows by symmetry.
\end{proof}

\smallskip

\begin{remark}
Recall that the log canonical threshold of $(X, V(\ua))$ is given by
$$c(X, V(\ua))=\sup\{\alpha\in\QQ_+\mid\cI(X,\alpha\cdot\ua)=\cO_X).$$
Since $\cI(X, \alpha\cdot\ua)=\cO_X$ if and only if
$\cI(X, c\cdot\ua)\cdot\cO_{V(\ua)}=\cO_{V(\ua)}$, it follows from
Proposition~\ref{invariance} that $\dim\,X-c(X, V(\ua))$ is independent
on the embedding $V(\ua)\hookrightarrow X$ (see \cite{mustata}
for an intrinsic expression for this difference). In fact, we show below
 that this is the case with all the jumping numbers of the multiplier
ideals of $\ua$.
\end{remark}

\medskip

If $\ua\subseteq\cO_X$ is a nonzero sheaf of ideals, then we say that
$\alpha\in\QQ$ is a jumping number of $\ua$ if
$$\cI(X, \alpha'\cdot\ua)\neq\cI(X, \alpha\cdot\ua),$$
for every $\alpha'<\alpha$. 

\begin{proposition}\label{jump}
With the notation in Proposition~\ref{invariance}, if there is an isomorphism
$\phi\,:\,Y_1\longrightarrow Y_2$, then for every $\gamma\in\QQ$,
$\gamma+\dim\,X_1$ is a jumping number for $\ua_1$ if and only if
$\gamma+\dim\,X_2$ is a jumping number for $\ua_2$.
\end{proposition}

\begin{proof}
The argument has the same flavor as the one used in
 Proposition~\ref{invariance}, so that we just sketch it briefly. 
Again. it is clear that we may assume that $\dim\,X_1=\dim\,X_2=n$.
The main point is to use Proposition~\ref{approximation} to show
that $\alpha$ is not a jumping number for a sheaf of ideals $\ua$
on a variety $X$  if and only if there is
$\epsilon\in\QQ_+^*$ such that for every $x\in X$ and every 
$p$, $q\geq 1$, we have
\begin{equation}\label{reduction}
\cI\left(X, (\alpha-\epsilon)\cdot (\ua+\um_x^p)\right)\subseteq
\cI\left(X, (\alpha+\epsilon)\cdot(\ua+\um_x^q)\right)
+\um_x^{[p\epsilon]-n+1}.
\end{equation}
Indeed, suppose first that $\alpha$ is not a jumping number.
Then there is $\epsilon\in\QQ_+^*$ such that
$$\cI\left(X, (\alpha-2\epsilon)\cdot\ua\right)=\cI\left(X,
(\alpha+\epsilon)\cdot\ua\right).$$ 
Proposition~\ref{approximation} gives
$$\cI\left(X,(\alpha-\epsilon)\cdot (\ua+\um_x^p)\right)
\subseteq\cI\left(X, (\alpha-2\epsilon)\cdot\ua\right)
+\um_x^{[p\epsilon]-n+1}=$$
$$\cI\left(X, (\alpha+\epsilon)\cdot\ua\right)+\um_x^{[p\epsilon]-n+1}
\subseteq\cI\left(X, (\alpha+\epsilon)\cdot (\ua+\um_x^q)\right)
+\um_x^{[p\epsilon]-n+1},$$
which gives the inclusion in equation~(\ref{reduction}). 

Conversely, if we have equation~(\ref{reduction}), then
Proposition~\ref{approximation} gives
$$\cI\left(X, (\alpha-\epsilon)\cdot(\ua+\um_x^p)\right)\subseteq
\cI(X, \alpha\cdot\ua)+\um_x^{[q\epsilon]-n+1}+\um_x^{[p\epsilon]-n+1},$$
for all $q\geq 1$. Krull's Intersection theorem implies  that
$$\cI\left(X, (\alpha-\epsilon)\cdot(\ua+\um_x^p)\right)\cdot\cO_{X,x}
\subseteq \cI(X, \alpha\cdot\ua)\cdot\cO_{X,x}
+\um_x^{[p\epsilon]-n+1}\cdot\cO_{X,x}.$$
Since the left hand side of the above inclusion contains
$\cI\left(X,(\alpha-\epsilon)\cdot\ua\right)\cdot\cO_{X,x}$,
one more application of Krull's Intersection theorem shows that
$\alpha$ is not a jumping number.

\smallskip 

Since equation~(\ref{reduction}) is a statement
 about
zero dimensional subschemes, the proof can be concluded with an argument
which paralels the one in
Proposition~\ref{invariance}.
\end{proof}

%References
\providecommand{\bysame}{\leavevmode \hbox \o3em
{\hrulefill}\thinspace}

\end{document}